\documentclass[12pt,a4paper,reqno]{article}

\usepackage[T1]{fontenc}
\usepackage[english]{babel}
\usepackage{amsmath}
\input xy \xyoption{all}

\newcounter{itm}
\newcounter{env}[section]
\renewcommand\theenv{\thesection.\arabic{env}}

\makeatletter

\renewcommand\l@section{\@dottedtocline{1}{0em}{1.5em}}

\def\section{%
  \if@noskipsec\leavevmode\fi
  \par
  \setlength\@tempskipa{-3.5ex\@plus-1ex\@minus-.2ex}
  \@afterindenttrue
  \ifdim\@tempskipa<\z@
    \@tempskipa-\@tempskipa\@afterindentfalse
  \fi
  \if@nobreak
    \everypar{}
  \else
    \addpenalty\@secpenalty\addvspace\@tempskipa
  \fi
  \@ifstar
    {\unnumberedsection}
    {\@dblarg{\numberedsection}}}

\def\numberedsection[#1]#2{%
  \refstepcounter{section}
  \protected@edef\@svsec{\@seccntformat{section}\relax}
  \setlength\@tempskipa{2.3ex\@plus.2ex}
  \ifdim\@tempskipa>\z@
    \begingroup
    \normalfont\Large\bfseries{\@hangfrom{\hskip\z@\relax\S\@svsec}\interlinepenalty\@M#2\@@par}
    \endgroup
  \else
    \def\@svsechd{\normalfont\Large\bfseries{\hskip\z@\relax\S\@svsec#2}}
  \fi
  \csname sectionmark\endcsname{#1}
  \addcontentsline{toc}{section}{\S\protect\numberline{\csname thesection\endcsname}#1}
  \@xsect{2.3ex\@plus.2ex}}

\def\unnumberedsection#1{%
  \setlength\@tempskipa{2.3ex\@plus.2ex}
  \ifdim\@tempskipa>\z@
    \begingroup
    \normalfont\Large\bfseries{\@hangfrom{\hskip\z@}\interlinepenalty\@M#1\@@par}
    \endgroup
  \else
    \def\@svsechd{\normalfont\Large\bfseries{\hskip\z@\relax#1}}
  \fi
  \@xsect{2.3ex\@plus.2ex}}

\def\subsection{\@ifstar\sottosezione\sottosezione}
\def\sottosezione#1{\@startsection{subparagraph}{2}{\z@}{-.7\baselineskip}{.3\baselineskip}{\bfseries}{#1}}

\renewcommand\appendix{%
  \par
  \setcounter{section}{0}
  \gdef\thesection{\@Alph\c@section}
  \renewcommand\numberedsection{\appendice}}

\def\appendice[#1]#2{%
  \refstepcounter{section}
  \protected@edef\@svsec{\thesection:~~}
  \setlength\@tempskipa{2.3ex\@plus.2ex}
  \ifdim\@tempskipa>\z@
    \begingroup
    \normalfont\Large\bfseries{\@hangfrom{\hskip\z@}Appendix~\@svsec\interlinepenalty\@M#2\@@par}
    \endgroup
  \else
    \def\@svsechd{\normalfont\Large\bfseries{\hskip\z@\relax Appendix~\@svsec#2}}
  \fi
  \csname sectionmark\endcsname{#1}
  \addcontentsline{toc}{section}{Appendix~\protect\numberline{\csname thesection\endcsname\relax:} #1}
  \@xsect{2.3ex\@plus.2ex}}

\makeatother

\newif\ifempty
\def\testifempty#1{\if#1\empty\relax\relax\emptytrue\else\emptyfalse\fi}

\providecommand\qedsymbol{\textit{q.e.d.}}
\providecommand\sdsymbol{\ensuremath{\blacksquare}}
\newcommand\mathqed{\quad\hbox{\qedsymbol}}
\newcommand\mathsd{\quad\hbox{\sdsymbol}}

\DeclareRobustCommand\qed{\ifmmode\mathqed\else\leavevmode\unskip\penalty9999\hbox{}\nobreak\hfill\quad\hbox{\qedsymbol}\fi}

\DeclareRobustCommand\sd{\ifmmode\mathsd\else\leavevmode\unskip\penalty9999\hbox{}\nobreak\hfill\quad\hbox{\sdsymbol}\fi}

\newenvironment{paragrafo}[1][]{\refstepcounter{env}\begin{list}{}{\setlength\itemindent{0pt}\setlength\labelsep\parindent\setlength\labelwidth{0pt}\setlength\leftmargin{0pt}\setlength\listparindent\parindent\setlength\parsep\parskip\setlength\partopsep{0pt}}\item\textbf{\theenv~~}\testifempty{#1}\ifempty\else\textbf{#1\ \ }\fi}{\end{list}}

\newenvironment{enunciato}[2][]{\refstepcounter{env}\begin{list}{}{\setlength\itemindent{0pt}\setlength\labelsep{0pt}\setlength\labelwidth{0pt}\setlength\leftmargin\parindent\setlength\listparindent\parindent\setlength\parsep\parskip\setlength\partopsep{0pt}}\item\textbf{\theenv~~}\testifempty{#2}\ifempty\else\textbf{#2~~}\fi\testifempty{#1}\ifempty\else(#1)~~\fi\slshape}{\end{list}}

\newenvironment{lemma}[1][]{\begin{enunciato}[#1]{Lemma}}{\end{enunciato}}
\newenvironment{proposizione}[1][]{\begin{enunciato}[#1]{Proposition}}{\end{enunciato}}

\newenvironment{enunciato*}[2][]{\begin{list}{}{\setlength\itemindent{0pt}\setlength\labelsep{0pt}\setlength\labelwidth{0pt}\setlength\leftmargin\parindent\setlength\listparindent\parindent\setlength\parsep\parskip\setlength\partopsep{0pt}}\item\testifempty{#2}\ifempty\else\textbf{#2~~}\fi\testifempty{#1}\ifempty\else(#1)~~\fi\slshape}{\end{list}}

\newenvironment{proof}[1][]{\begin{list}{}{\setlength\itemindent\parindent\setlength\labelsep{0pt}\setlength\labelwidth{0pt}\setlength\leftmargin{0pt}\setlength\listparindent\parindent\setlength\parsep\parskip\setlength\partopsep{0pt}}\item\testifempty{#1}\ifempty\emph{Proof.}~~\else\emph{#1.}\ \fi}{\qed\end{list}}

\newenvironment{elenco}{\begin{list}{(\roman{itm})}{\setlength\itemindent{0pt}\setlength\labelsep{0.5em}\setlength\labelwidth\leftmargin\addtolength\labelwidth{-\labelsep}\setlength\listparindent{.5\parindent}\setlength\parsep\parskip\setlength\itemsep\medskipamount\setlength\partopsep{0pt}\usecounter{itm}}}{\end{list}}

\usepackage{amssymb}
\usepackage{mathrsfs}
\usepackage{bbold}

\newcommand\Tau{{\mathrm T}}

\newcommand\inciso[1]{\nobreakdash---\hspace{0pt}{#1}\nobreakdash---\hspace{0pt}}

\newcommand\textsum[3][]{\ensuremath{\overset{#1}{\underset{#2}{\textstyle\sum}}\,#3}}

\newcommand\textcup[3][]{\ensuremath{\overset{#1}{\underset{#2}{\textstyle\cup}}\,#3}}

\newcommand\nN{\ensuremath{\mathbb{N}}}

\newcommand\nR{\ensuremath{\mathbb{R}}}
\newcommand\nC{\ensuremath{\mathbb{C}}}

\newcommand\fro\leftarrow
\newcommand\infro\hookleftarrow
\newcommand\into\hookrightarrow

\newcommand\isoto{\stackrel\thicksim\to}
\newcommand\longfro\longleftarrow
\newcommand\longto\longrightarrow
\newcommand\onfro\twoheadleftarrow
\newcommand\onto\twoheadrightarrow

\newcommand\xto[1]{\xrightarrow{#1}}

\newcommand\can\cong
\newcommand\iso\approx

\newcommand\id{\mathit{id}}
\newcommand\Id{\mathit{Id}}

\newcommand\algdim[2][]{\testifempty{\ensuremath{#1}}\ifempty\ensuremath{\mathrm{dim}\,#2}\else\ensuremath{\mathrm{dim}_{#1}#2}\fi}

\newcommand\kernel[1]{\ensuremath{\mathrm{Ker}\,#1}}

\newcommand\transposed[1]{\ensuremath{{^{\mathsf t}#1}}}

\newcommand\End{\mathrm{End}}
\newcommand\Aut{\mathrm{Aut}}

\newcommand\GL{\mathit{GL}}
\newcommand\OG{\mathit{O}}
\newcommand\SO{\mathit{SO}}

\newcommand\so{\mathfrak{so}}

\newcommand\supp[2][]{\testifempty{\ensuremath{#1}}\ifempty\ensuremath{\mathrm{supp}\,{#2}}\else\ensuremath{\mathrm{supp}_{#1}{#2}}\fi}

\newcommand\norm[2][]{\testifempty{\ensuremath{#1}}\ifempty\ensuremath{\left|{#2}\right|}\else\ensuremath{\left\|{#2}\right\|_{#1}}\fi}

\newcommand\supnorm[1]{\ensuremath{\left\|#1\right\|}}

\newcommand\C{\mathit C}

\newcommand\T[2][]{\testifempty{\ensuremath{#1}}\ifempty\ensuremath{\mathit T\,#2}\else\ensuremath{\mathit T_{#1}#2}\fi}

\newcommand\aeq\Leftrightarrow
\newcommand\seq\Rightarrow

\newcommand\s[1][]{\testifempty{\ensuremath{#1}}\ifempty\ensuremath{\boldsymbol s}\else\ensuremath{\boldsymbol s\mspace{.7mu}#1}\fi}

\renewcommand\t[1][]{\testifempty{\ensuremath{#1}}\ifempty\ensuremath{\boldsymbol t}\else\ensuremath{\boldsymbol t\mspace{.7mu}#1}\fi}

\newcommand\pprin[1]{#1^\mathrm{pr}}
\newcommand\psing[1]{#1^!}

\newcommand\e[1]{\ensuremath{\boldsymbol e_{#1}}}

\begin{document}

\title{A~note on~the~global~structure of~proper~Lie~groupoids in~low~codimensions\thanks{Research partly supported by Fon\-da\-zio\-ne "Ing.~Al\-do Gi\-ni"}}
\author{Giorgio Trentinaglia}
\date{\itshape\footnotesize Courant Research Centre ``Higher Order Structures''\\University of G\"ottingen, Germany\thanks{Supported by Deutsche Forschungs\-gemein\-schaft (DFG; German Research Foundation) in the context of the strategic excellence program of the University of G\"ottingen}}
\maketitle

\begin{abstract}
\noindent We observe that any connected proper Lie groupoid whose orbits have codimension at most two admits a \textit{globally effective} representation on a smooth vector bundle, i.e.\ one whose kernel consists only of ineffective arrows. As an application, we deduce that any such groupoid can up to Morita equivalence be presented as an extension of some action groupoid $G \ltimes X$ with $G$ compact by some bundle of compact Lie groups.
\end{abstract}

\section*{Motivation}
A typical problem in the theory of Lie groupoids is to decide whether a Lie groupoid can be reduced, up to Morita equivalence, to a certain standard form\nobreakdash---\hspace{0pt}for instance to a translation groupoid $G \ltimes X$ for some smooth action of a compact Lie group $G$ on a manifold $X$. Problems of this sort are usually called \textit{presentation problems} in the literature \cite{HeMe04}. The importance of this kind of question is self-evident, as, conventionally, the ``{standard form}'' groupoids one has in mind are always familiar and well studied mathematical objects.

For proper Lie groupoids, being Morita equivalent to a translation groupoid for some compact Lie group action is, in general, a severe restriction. In fact, Lie groupoids of fairly simple type like bundles of compact Lie groups over a connected, compact manifold may fail to have that property; a standard example can be found in \cite[\S2.10]{2008b}. However, when one is only interested in the effective geometry of proper Lie groupoids, i.e.\ in the study of the global properties of the canonical action of these groupoids on their own base manifolds, it may be reasonable to see whether less stringent presentability problems are solvable in greater generality. For instance, one might ask whether, taken for granted the admissibility of factoring out any \textit{ineffective} subbundle of a proper Lie groupoid, it would be possible in general to get a quotient Morita equivalent to a translation groupoid for some compact Lie group action (an isotropic arrow of a Lie groupoid is said to be \textit{ineffective} if its infinitesimal action on the normal space to the corresponding orbit is trivial; the precise definition is given in \S\ref{sez1} below). More precisely, consider the following problem: given a proper Lie groupoid $\mathcal G$, find a short exact sequence of Lie groupoid homomorphisms
\begin{equation}\label{equ:pre-seq}
1\; \to\; \mathcal K\; \into\; \mathcal G\; \stackrel\pi\onto\; \mathcal G'\; \to 1
\end{equation}
where $\pi$ is a surjective submersion identical on the base manifolds, so that (a)~$\mathcal K$ is an ineffective subbundle of $\mathcal G$, and (b)~$\mathcal G'$ is Morita equivalent to an action groupoid $G \ltimes X$ for some action of a compact Lie group $G$ on a manifold $X$. (The first property is actually implied by the second one; compare \cite[Lemma 4.6 and Prop.~4.1]{2008b}.) Of course one should like the quotient $\mathcal G'$ to differ as little as possible from the original groupoid $\mathcal G$, so one might refine the problem by requiring of the epimorphism $\mathcal G \onto \mathcal G'$ the following universal property:
\begin{equation}\label{equ:pre-univ-prop}
\begin{split}
\xymatrix@C=70pt{\mathcal G\ar[dr]_{\forall\:\varphi}\ar@{->>}[r]^-\pi & \mathcal G'\ar@{-->}[d]^{\:\exists\:!\:\varphi'} \\ & \mathcal H}
\end{split}
\end{equation}
for each Lie groupoid homomorphism $\varphi$ into a Lie groupoid $\mathcal H$ that is Morita equivalent to a translation groupoid $H \ltimes Y$ for some compact Lie group action $H \circlearrowright Y$. However, one can show that if at least one sequence (\ref{equ:pre-seq}) can be constructed then, when the orbit space of the groupoid has finitely many components, there is also some sequence for which the universal property (\ref{equ:pre-univ-prop}) holds. In the present paper we intend to point out a simple sufficient condition for the solvability of the above presentation problem (under some reasonable hypothesis e.g.\ connectedness of the orbit space), namely, that the singular foliation of the base manifold of $\mathcal G$ given by the $\mathcal G$-orbits be of codimension at most two. It is our hope that the study of presentation problems like (\ref{equ:pre-seq}) will pave the way to classification results for proper Lie groupoids in the spirit of \cite{Moe03}.

A standard result in the representation theory of Lie groupoids states that a necessary and sufficient condition in order that a proper Lie groupoid $\mathcal G$ may be Morita equivalent to a translation groupoid $G \ltimes X$ with $G$ compact is the existence of a faithful representation $\varrho: \mathcal G \to \GL(E)$ on some vector bundle $E$ of constant rank over the base of $\mathcal G$ (the reader can find a proof in \cite[Section 5]{2008}). Similarly, in order to find a solution to the presentation problem (\ref{equ:pre-seq}) for $\mathcal G$, it will be enough to construct an \textit{effective} representation of $\mathcal G$ on some vector bundle of constant rank\nobreakdash---a representation being \textit{effective} when its kernel is contained in the ineffective sub\-bundle of $\mathcal G$. Other problems related to the existence of specific types of representation for a generic proper Lie groupoid originate from the duality theory developed in \cite{2008b}. Let us call a proper Lie groupoid $\mathcal G$ \textit{parareflexive} when for each base point $x$ there is a representation whose kernel at $x$ is ineffective; also, call $\mathcal G$ \textit{reflexive} when the analogous condition obtained by replacing the word `ineffective' with the word `trivial' holds. We have at present no counterexamples to any of the following statements:

\smallskip
(S1) {\em Every proper Lie groupoid is parareflexive.}

\smallskip
(S2) {\em Every proper Lie groupoid admits an effective representation on a vector bundle of constant rank.} Equivalently, every such groupoid can up to Morita equivalence be presented as an extension of some action groupoid $G \ltimes X$ with $G$ compact by some bundle of compact Lie groups.

\smallskip
(S3) {\em An arbitrary connected proper Lie groupoid is reflexive if, and only if, it admits a faithful representation on a vector bundle of constant rank} (equivalently: if, and only if, it is Morita equivalent to a translation groupoid for some compact Lie group action).

\smallskip
In fact, (S3) is an immediate corollary of (S2) in the connected case. If counterexamples to (S1) exist then\inciso{as we are going to show in the present paper}they ought to live in base dimensions greater or equal to three. Even in codimension two, we do not know whether (S2) holds for a groupoid with infinitely many connected components.

\section{Preliminaries}\label{sez1}
We assume that the reader has some familiarity with the general theory of Lie groupoids, at the level of let us say Chapter 5 of \cite{MoeMrc03}.

We start by recalling the notion of ineffectiveness for isotropic arrows of a Lie groupoid. Let $\mathcal G$ be an arbitrary Lie groupoid, $X$ its base manifold, and $x \in X$ a base point. For each arrow $g \in \mathcal G_x$ in the isotropy group of $\mathcal G$ at $x$ one obtains a well-defined linear automorphism of the space ${N^{\mathcal G}}_x = {\T[x]{X} / \T[x]{O_x}}$ of all tangent vectors perpendicular to the orbit $O_x = {\mathcal G\cdot x}$ by first choosing an arbitrary local bisection $\sigma: U \into \mathcal G$, ${\t \circ \sigma}: U \isoto U'$ with $\sigma(x) = g$ and then taking the quotient linear map induced on ${N^{\mathcal G}}_x$ by the tangent map $\T[x]{(\t \circ \sigma)}$. In fact, one obtains a continuous representation of the Lie group $\mathcal G_x$ on the vector space ${N^{\mathcal G}}_x$, which we shall denote by $\mu_x$. When $g$ belongs to $\kernel{\mu_x}$, one calls $g$ \textit{ineffective.}

Let $\varrho$ be a representation of $\mathcal G$ on a smooth vector bundle $E$ of constant rank over $X$, that is, a homomorphism of Lie groupoids $\varrho: \mathcal G \to \GL(E)$ inducing the identity on $X$. We shall say that $\varrho$ is \textit{effective at $x$} when the kernel of the induced isotropy representation $\varrho_x: \mathcal G_x \to \GL(E_x)$ is ineffective, i.e., contained in $\kernel{\mu_x}$. We shall call $\varrho$ \textit{effective} when it is effective at all $x$.

It is our purpose in this note to prove that every connected proper Lie groupoid whose orbits have codimension at most two admits an effective representation (recall that a Lie groupoid is said to be \textit{connected} when its orbit space is a connected topological space). We will give a direct and rather ad hoc construction which on the one hand has the advantage of showing that in some cases one can actually find an effective representation $\varrho$ such that $\kernel{\varrho_x}$ is not larger than the connected component of the identity in the group $\mathcal G_x$ for all $x$, but on the other hand has the drawback that it does not carry over to the situation where $\mathcal G$ has infinitely many connected components. However, in view of the following example, the connectedness hypothesis appears justified to us at least in relation to the study of the presentation problem described at the beginning:

\sloppy
\begin{paragrafo}[Example]\label{exm:symmgrps}
Consider the bundle of finite groups over the discrete manifold $\nN = \{0, 1, 2, \ldots\}$ whose isotropy group at $n$ is the symmetric group $\mathfrak S_n$. For $n>6$, the only irreducible representations of $\mathfrak S_n$ of rank less than $n-1$ are the alternating and the trivial one; compare \cite[Exercise 4.14, p.~50]{FuHa91}. Hence, although all representations of this groupoid are effective, none of them can be ``universal'' in the sense expressed by (\ref{equ:pre-univ-prop}) above.
\end{paragrafo}

\fussy
Let $\mathcal G \rightrightarrows M$ be a (connected) proper Lie groupoid. We shall let $\pprin{M}$ denote the set of all base points $x \in M$ such that the whole isotropy group of $\mathcal G$ at $x$ is ineffective, i.e., such that $\kernel{\mu_x} = \mathcal G_x$. From the explicit description of $\mu_x$ given above it follows immediately that $\pprin{M}$ is a $\mathcal G$-invariant subset of $M$. Moreover, by the properness of $\mathcal G$, $\pprin{M}$ must be open; this is an obvious consequence of the local linearizability theorem for proper Lie groupoids \cite[Theorem 2.3]{Zu06}. We shall denote by $\pprin{\mathcal G} \rightrightarrows \pprin{M}$ the Lie groupoid induced on $\pprin{M}$ by restriction, and call this the \textit{principal part} of $\mathcal G$. Our terminology here is in agreement with \mbox{Bredon's} \cite[\S IV.3 and especially Theorem IV.3.2 (iii)]{Bre72}.

Incidentally, we observe that $\pprin{M}$ must be dense in $M$. This follows immediately from the local linearizability theorem and \cite[Theorem IV.3.1]{Bre72}; for any linear slice $i: V \into M$, $\pprin{M} \cap i(V)$ must be relatively dense in $i(V)$.

Of course $\pprin{\mathcal G} \rightrightarrows \pprin{M}$ is a regular proper Lie groupoid. Hence the isotropy bundle $\mathit I(\pprin{\mathcal G})$ is an embedded submanifold of $\pprin{\mathcal G}$ and, in fact, a locally trivial bundle of compact Lie groups over $\pprin{M}$ \cite{Moe03}. Note that by the local linearizability theorem and by \cite[Theorem IV.3.1]{Bre72} connectedness of $\mathcal G$ implies connectedness of $\pprin{\mathcal G}$, so the fibres of $\mathit I(\pprin{\mathcal G}) \to \pprin{M}$ are all isomorphic (as Lie groups).

\begin{paragrafo}[Component representation]\label{par:comprepr}
Let $\mathcal K \xto{\s=\t} M$ be any bundle of compact Lie groups. Suppose that all the isotropy groups of $\mathcal K$ are isomorphic to a fixed compact Lie group $K$. Let $\pi_0(K)$, and, for each $x \in X$, $\pi_0(\mathcal K_x)$, denote the finite group obtained by factoring out the identity connected component. Then there is a canonical representation of $\mathcal K$ on a vector bundle $R \equiv R^{\mathcal K}$ over $X$, of rank equal to the order $N$ of $\pi_0(K)$, to be called the \textit{component representation} of $\mathcal K$ and to be denoted by $\varrho \equiv \varrho^{\mathcal K}$, defined as follows.

For each base point $x \in X$, let
\begin{equation}\label{comprepr1}
R_x\, :=\, \C^0 \bigl( \pi_0(\mathcal K_x), \nC \bigr)\, \iso\, \nC^N
\end{equation}
be the vector space of all complex functions on $\pi_0(\mathcal K_x)$. The local triviality of $\mathcal K \to X$ yields an evident smooth vector bundle structure on $R = {\coprod R_x} \to X$. Let $k \in \mathcal K_x$ act on $R_x$ by right translation:
\begin{equation}\label{comprepr2}
\varrho(k)(f)\, :=\, {f \circ \pi_0(\tau^k)}
\end{equation}
for all $f: \pi_0(\mathcal K_x) \to \nC$, where $\pi_0(\tau^k)$ denotes the permutation of $\pi_0(\mathcal K_x)$ induced by the right translation diffeomorphism $h \mapsto hk$ of $\mathcal K_x$ onto itself.
\end{paragrafo}

\begin{lemma}\label{lem:eta}
Let $\mathcal H \rightrightarrows B$ be a connected, principal, proper Lie groupoid. Suppose that the orbit foliation of the base manifold $B$ has codimension one. Then there exists a representation $\eta: \mathcal H \to \GL(E)$ which on the isotropy bundle $\mathit I(\mathcal H) \to B$ induces, locally up to isomorphism, the component representation $\varrho^{\mathit I(\mathcal H)}$.
\end{lemma}\begin{proof}
Since $\mathcal H \rightrightarrows B$ is a principal Lie groupoid, the orbit space $X = B/\mathcal H$ endowed with the evident functional structure \cite[\S VI.1]{Bre72} is a smooth manifold. In fact, any principal linear slice $S \subset B$ determines a $\C^\infty$ parametrization of the open subset $\phi(S)$ of $X$ via the quotient projection $\phi: B \to X$, hence by our assumptions $X$ is one-dimensional. Moreover, since $\mathcal H$ is proper, $X$ must be Hausdorff.

Select a sequence $\{S_i\}_{i=1,2,\ldots}$ of principal linear slices ($S_i \simeq \nR$ and $\mathcal H|_{S_i} \simeq H \times \nR$ for some compact Lie group $H$) with $B = \textcup[\infty]{i=1}{\mathcal H\cdot S_i}$. For each $k$ put $B_k := \textcup[k]{i=1}{\mathcal H\cdot S_i}$. By rearranging the sequence if necessary, we may assume $S_{k+1} \cap B_k \neq \varnothing$ for all $k$ (here we use the connectedness of $\mathcal H$).

Next, select an invading sequence of open subsets $\ldots \subset V_p \subset \overline V_p \subset V_{p+1} \subset \ldots \subset B$ with $\phi(\overline V_p)$ compact. It is no loss of generality to assume $V_p$ to be invariant. For eack $k$ let $p(k)$ be the greatest $p \leqq k$ such that $\overline V_p \subset B_k$.

(Induction step.) Suppose a representation $\eta_k: \mathcal H_k \to \GL(E_k)$ as in the statement of the lemma has been constructed for $\mathcal H_k \equiv \mathcal H|_{B_k}$. Put $S \equiv S_{k+1}$ and $V \equiv V_{p(k)}$. By the initial remarks, the intersection $S \cap B_k$ has the form $(-\infty,a) \cup (b,\infty)$ or $(a,b)$ in any parametrization $S \simeq \nR$. It is then clear that we can find open subsets $\Sigma \subset S$ and $U \subset B_k$ with $\overline V \subset U$, $\overline V\cap \Sigma = \varnothing$, $U\cup {(\mathcal H\cdot \Sigma)} = B_k\cup {(\mathcal H\cdot S)}$ such that there is an isomorphism of representations on $U\cap \Sigma$ between $E_k$ and $R$. By a standard Morita equivalence argument, we obtain a representation $\eta_{k+1}: \mathcal H_{k+1} \to \GL(E_{k+1})$ which still satisfies the induction hypothesis and moreover is globally isomorphic to $\eta_k$ over $V$.

Finally, define the representation $\eta$ on $E$ to be the inductive limit of the partial representations $\eta_k|V_{p(k)}$ on $E_k|V_{p(k)}$.
\end{proof}

{\em Remark.} The above proof makes essential use of the codimension~= 1 hypothesis and of the fact that the component representation is defined intrinsically in terms of the groupoid structure.

\medskip

Let $\mathcal G \rightrightarrows M$ be a connected, proper Lie groupoid. Assume that the orbits of $\mathcal G$ have codimension at most two, in other words,
\begin{equation}\label{cond:codim2}
{\sup_{m\in M} \algdim{\bigl(\T[m]{M} / \T[m]{O_m}\bigr)}} \leqq 2\text.
\end{equation}
By the connectedness of $\mathcal G$, the dimension of the normal space $N_m = {\T[m]{M} / \T[m]{O_m}}$ is the same for all $m \in \pprin{M}$. We shall call the common dimension of these spaces the \textit{principal codimension} of $\mathcal G$. The latter equals the dimension of the manifold $\pprin{X} := {\pprin{M} / \pprin{\mathcal G}}$. (Observe that $\pprin{X}$ is a dense open subspace of the orbit space $X = M/\mathcal G$.) We will now indicate how to obtain an effective representation for $\mathcal G$ in different ways according to whether the principal codimension of the groupoid equals zero, one or two.

\smallskip
(Principal codimension zero.) In this case the singular set $\psing{M} := {M \setminus \pprin{M}}$ is empty and $\mathcal G$ is transitive. Hence there actually is a {\em faithful} representation, by Morita equivalence to a compact Lie group.

\smallskip
(Principal codimension one.) Let $m \in \psing{M}$ be a singular base point, and consider any linear slice $\nR^d \simeq V \subset M$, $\mathcal G|_V \simeq G\ltimes \nR^d$ at $m$ for some let us say orthogonal action $G \to \GL(\nR^d)$ of the isotropy group $G \equiv \mathcal G_m$. One then has the following possibilities:
\begin{elenco}
\item $d=1$, and $G \to \OG(\nR) = \{\pm\Id\}$ is a nontrivial action with orbits $\{\pm t\}$ where $t \in \nR$;
\item $d=2$, and $G \to \OG(2) := \OG(\nR^2)$ is a \textit{spherical} action, i.e., one whose orbits are the circles $x^2+y^2 = r^2$ with $r \geqq 0$;
\end{elenco}
from this remark it follows in particular that $\psing{X} := {X \setminus \pprin{X}}$ is a discrete subset of the orbit space $X = {M/\mathcal G}$, that is, all its points are isolated in $X$. Hence we can assign each $x \in \psing{X}$ an open neighbourhood $\Omega_x$ so that $\Omega_x \cap \Omega_{x'} = \varnothing$ for $x \neq x' \in \psing{X}$. Let us for each $x \in \psing{X}$ fix some $m \equiv m_x \in \psing{M}$ with $x = \phi(m)$, where $\phi: M \to X$ denotes the projection onto the quotient, and some (orthogonal) linear slice $V \equiv V_x$ at $m$ with $\phi(V) \subset \Omega_x$.

We contend that for each orthogonal compact Lie group action $G \to \OG(\nR^d)$ as in (i) or (ii) one can find an effective representation $\Phi: G\ltimes \nR^d \to \GL(\underline\nC^{2N})$ (where $\underline\nC^{2N} = \nR^d \times \nC^{2N}$) such that its restriction to the isotropy bundle over $\pprin{(\nR^d)} = {\nR^d \setminus 0}$ is locally isomorphic to twice the component representation. This is clear in the case (i) (compare Remark \ref{rmk:spher-rk1} below) and shall be proved in the next section for actions of type (ii). Taking for granted the existence of such representations, fix one for each linear slice $V_x$ and name it $\Phi_x$. Moreover choose a representation $\eta: \pprin{\mathcal G} \to \GL(E)$ as in the statement of Lemma \ref{lem:eta}. Then there are invariant open subsets $U \subset \pprin{M}$ and $B_x \subset V_x$ such that $U \cup {\bigcup \bigl\{U_x: x \in \psing X\bigr\}} = M$, where $U_x \equiv {\mathcal G\cdot B_x}$, and $U\cap U_x = {\mathcal G\cdot \Sigma}$ for some $\Sigma \equiv \Sigma_x \subset V_x$ of the form $(a,b) \times \nR^{d-1} \subset \nR^d$ ($a<b$ positive) with $\Phi_x|_\Sigma \iso R^{\mathcal G|_\Sigma} \oplus R^{\mathcal G|_\Sigma} \iso \eta|_\Sigma \oplus \eta|_\Sigma$ as representations of $\mathcal G|_\Sigma = \mathit I(\mathcal G|_\Sigma)$. Now, the usual Morita equivalence techniques allow one to glue together these representations into a global representation which has to be isomorphic to $\eta \oplus \eta$ on $U$ and to $\Phi_x$ on $B_x$ and therefore is effective.

\smallskip
(Principal codimension two.) In this case the groupoid is {\em regular,} hence the result follows immediately from \cite[Corollary 4.12]{2008b}.

\begin{paragrafo}[Remark: spherical orthogonal actions of rank one]\label{rmk:spher-rk1}
Let $\mu: G \to \OG(1) = \{\pm I\}$ be a rank-one orthogonal linear action of a compact Lie group $G$ with spherical orbits $\{\pm t\}$. Put $K \equiv \kernel{\mu}$. Clearly $K = G_t$ for all $t \neq 0$, where $G_t$ denotes the stabilizer at $t$. Also, the connected component of the identity in $G_t$ and in $G$ are the same; in symbols, $(G_t)^{(e)} = G^{(e)}$. Hence it will be no loss of generality to assume that $G$ is discrete.

Now let $\varrho^G: G \into \GL(R^G)$ be the (right) regular representation of the (finite) group $G$. One has $R^G = \C^0(G,\nC) = \C^0(K,\nC) \oplus \C^0({G \setminus K},\nC) \iso R^K \oplus R^K$ equivariantly as $K$-modules. The obvious extension $\Phi: G\ltimes \nR \to \GL \bigl(\nR \times R^G\bigr)$ of $\varrho^G$ is then an effective (in fact, faithful) representation whose restriction to the isotropy subbundle over ${\nR \setminus 0}$ is isomorphic to twice the component representation.
\end{paragrafo}

\section[Spherical orthogonal actions of rank two]{Spherical orthogonal actions on $\nR^2$}\label{sez2}
Let $G$ be a compact Lie group acting orthogonally and spherically on $\nR^2$, in other words, let a continuous homomorphism $\mu: G \to \OG(2) = \OG(\nR^2)$ be given such that the orbits of the corresponding action $G \circlearrowright \nR^2$ are the circles $x^2+y^2 = r^2$, $r\geqq0$. Our goal in this section is to construct an effective representation $\Phi: G\ltimes \nR^2 \to \GL(\underline\nC^{2N})$ such that its restriction to the isotropy bundle over ${\nR^2 \setminus 0}$ is locally isomorphic to the direct sum of two copies of the component representation. Clearly, it will be enough to find a continuous homomorphism of groups $\varphi: G \to \GL(\nC^{2N})$ such that $\kernel{\varphi} \subset K$ and such that $\varphi|_{G_1} \iso \varrho^{G_1} \oplus \varrho^{G_1}$ where $G_1$ denotes the stabilizer subgroup at $(1,0) \in \nR^2$, for then $\Phi$ defined by
$$%
{\Phi(g,z) \cdot (z,\boldsymbol v)} := \bigl(\mu(g)(z),\varphi(g)\boldsymbol v\bigr)
$$%
will have the desired properties.

We start with some remarks. Let $G^{(e)}$ denote the identity component of $G$ and $\mu^{(e)}: G^{(e)} \to \SO(2)$ the restriction of $\mu$ to $G^{(e)}$. Our first remark is that $G^{(e)} \cap K = G^{(e)} \cap G_1$ where $K = \kernel{\mu}$ and $G_1$ denotes the stabilizer at $(1,0) \in \nR^2$. Indeed, for every $x \in G^{(e)}$ the matrix $\mu^{(e)}(x) \in \SO(2)$ has no nonzero fixed vectors unless $x \in K$. Next, we observe that the composite $G_1 \subset G \to G/G^{(e)}$ is a surjective homomorphism of groups. Indeed, if $U$ is a connected component of $G$ and $g_0 \in U$ then ${{g_0}^{-1} \cdot U} = G^{(e)}$, hence ${U \cdot (1,0)} = {g_0 \cdot G^{(e)} \cdot (1,0)} = {g_0 \cdot S^1} = S^1$ by sphericity and so there exists some $g \in U$ with ${g \cdot (1,0)} = (1,0)$, which tells us that $U \cap G_1 \neq \varnothing$. From the first remark it follows that $(G_1)^{(e)} \subset G^{(e)} \cap K$ and therefore that $(G_1)^{(e)} = \smash{\bigl(G^{(e)} \cap K\bigr)}^{(e)}$ is a normal subgroup of $G$ contained in $K$. Hence by factoring out $(G_1)^{(e)}$ we are reduced to the special case where the stabilizer subgroup $G_1$ is {\em finite.} This will be our assumption for the remainder of the section.

By sphericity, the homomorphism $\mu^{(e)}: G^{(e)} \to \SO(2) = S^1$ must be surjective. In particular, $\mu^{(e)}$ must be nontrivial, hence submersive. The kernel $\kernel{\mu^{(e)}} = G^{(e)} \cap K = G^{(e)} \cap G_1$ is finite. It follows that $\mu^{(e)}$ is a finite-sheeted covering of $S^1$ and that $G^{(e)}$ is one-dimensional. Let $\alpha: \nR \to G^{(e)}$ be the unique Lie group homomorphism such that ${\mu^{(e)} \circ \alpha} = \exp: \nR \to S^1$. Being surjective, $\alpha$ induces an isomorphism of Lie groups $S^1 \iso G^{(e)}$. Thus the subgroup $C \equiv G^{(e)} \cap K \subset G^{(e)}$ must be cyclic. Let us pick a generator:
\begin{equation}\label{def:cyclic-gen}
C = \langle c_0 \rangle = \bigl\{c_0, {c_0}^2, \ldots, {c_0}^q=e\bigr\}
\end{equation}
where $q = \norm{C}$ is the order. Let $H$ denote the kernel of the restriction $\mu|_{G_1}: G_1 \to \OG(2)$, and put $H' \equiv {G_1 \setminus H}$. We will now show that for each $g \in G_1$
\begin{equation}\label{equ:action-by-conj}
\begin{cases}
gxg^{-1} = x & \text{for all $x \in G^{(e)}$ if and only if $g \in H$}
\\
gxg^{-1} = x^{-1} & \text{for all $x \in G^{(e)}$ if and only if $g \in H'$}\text.
\end{cases}
\end{equation}
To begin with, $G^{(e)} \iso S^1$ implies $\Aut \bigl( G^{(e)} \bigr) = \{\id, \chi\}$ where $\chi$ denotes the inversion $x \mapsto x^{-1}$. Since $G^{(e)}$ is normal in $G$, each $g \in G_1$ acts on $G^{(e)}$ by conjugation and thus induces an automorphism $c_g \in \Aut \bigl( G^{(e)} \bigr)$, hence either $gxg^{-1} = x$ for all $x \in G^{(e)}$ or $gxg^{-1} = x^{-1}$ for all $x \in G^{(e)}$. Suppose first $g \in H$. By definition of $H$, $\mu(g) = \bigl( \begin{smallmatrix} 1 & 0 \\ 0 & 1 \end{smallmatrix} \bigr)$. Therefore, if $gxg^{-1} = x^{-1}$ $\forall x \in G^{(e)}$ then $\mu(x)^{-1} = \mu(x)$ $\forall x \in G^{(e)}$ and hence by connectedness $\mu(x) = \bigl( \begin{smallmatrix} 1 & 0 \\ 0 & 1 \end{smallmatrix} \bigr)$ $\forall x \in G^{(e)}$ which contradicts the finiteness of $G^{(e)} \cap K$. Suppose on the other hand that $g \in H'$. Then $\mu(g) = \bigl( \begin{smallmatrix} 1 & \phantom-0 \\ 0 & -1 \end{smallmatrix} \bigr)$. If one assumes that $gxg^{-1} = x$ $\forall x \in G^{(e)}$ then one gets a contradiction as before.

The induced representation $\mu|_{G_1}: G_1 \to \OG(2)$ embeds into the direct sum of two copies of the (real) regular representation $R \oplus R \equiv R^{G_1} \oplus R^{G_1}$ as the two-dimensional submodule
\begin{equation}\label{def:Tau}
\Tau := \Tau_1 \oplus \Tau_2 := \biggl\langle \textsum{g \in G_1}{\!\e{g}} \biggr\rangle \oplus \biggl\langle \textsum{h \in H}{\!\e{h}} - \textsum{h' \in H'}{\!\e{h'}} \biggr\rangle \subset R \oplus R
\end{equation}
where \e{g} denotes the standard basis vector given by the function with value one at $g \in G_1$ and zero everywhere else. Identify $R \equiv \C^0(G_1, \nR) = \mathit L^2(G_1, \nR)$ where $G_1$ has the probability Haar measure. The orthogonal complement of the submodule (\ref{def:Tau}) is then
\begin{multline}\label{def:Lambda}
\Lambda := \Lambda_1 \oplus \Lambda_2 := \biggl\{f \in \C^0(G_1)\biggl| \textsum{g \in G_1}{\!f(g)} = 0\biggr\}\; \oplus
\\
\biggl\{f \in \C^0(G_1)\biggl| \textsum{h \in H}{\!f(h)} = \textsum{h' \in H'}{\!f(h')}\biggr\} \subset R \oplus R
\end{multline}
where $R = \Tau_i \oplus \Lambda_i$ ($i = 1,2$). We now proceed to show how to extend the action of $G_1$ on $\Lambda$ to an action of the whole group. (In the end we shall define $\varphi$ as the complexified direct sum of $\mu$ and this extended representation.) We distinguish two cases.

\begin{paragrafo}[The odd-order case]\label{par:q=2p+1}
Suppose $q = \norm{C} = 2p+1$ for some $p \in \nN$. Let us put $N = \norm{G_1}$ (recall that $G_1$ is a finite group) and $n = N-1$.

Choose an arbitrary orthonormal basis for the $G_1$-submodule (\ref{def:Lambda}) $\Lambda \subset R \oplus R$, and let $\Pi_1: G_1 \to \OG\bigl( 2N-2 \bigr) = \OG(2n)$ be the corresponding orthogonal representation. Note that $N$ coincides with the rank of the $G_1$-module $R$. We contend that $\Pi_1$ can be extended to a continuous representation $\Pi: G \to \OG(2n)$. Of course, we may assume $n \geqq 1$.

To begin with, we remark that $\Pi_1(C) \subset \SO(2n)$. Indeed, on the one hand, from the inclusion $C = G^{(e)} \cap K \subset K$ we deduce $\det \mu(c) = \det(\id) = +1$ for all $c \in C$. On the other hand, from the identity (of $G_1$-modules) $\Tau \oplus \Lambda = R \oplus R$ we obtain
$$%
\det\mu(c)\, \det\Pi_1(c) =\, \det\begin{bmatrix} \mu(c) & 0 \\ 0 & \Pi_1(c) \end{bmatrix} =\, \det\begin{bmatrix} \varrho(c) & 0 \\ 0 & \varrho(c) \end{bmatrix} = \bigl[ \det\varrho(c) \bigr]^2 = +1
$$%
(because $\varrho: G_1 \to \OG(N)$ is orthogonal), so that, in fact, $\det\Pi_1(c) = +1$ for all $c \in C$.

Let $c_0$ be, as in (\ref{def:cyclic-gen}), the selected generator of the cyclic group $C$. Put $P_0 = \Pi_1(c_0) \in \SO(2n)$. Then ${P_0}^q = I_{2n}$ (identity $2n \times 2n$ matrix). Every element of $\SO(2n)$ is conjugated in $\SO(2n)$ to an element of the standard maximal torus $T(n)$ consisting of all block-diagonal matrices of the form $R(\theta_1, \dotsc, \theta_n) := \mathrm{diag} \bigl( R(\theta_1), \dotsc, R(\theta_n) \bigr)$, $\theta_1, \dotsc, \theta_n \in \nR$, where $R(\theta)$ denotes the $2 \times 2$ matrix $\bigl( \begin{smallmatrix} \cos\theta & -\sin\theta \\ \sin\theta & \phantom-\cos\theta \end{smallmatrix} \bigr) \in \SO(2)$. Compare \cite[Theorem (3.4)]{BtD95}. Thus, at the expense of replacing the representation $\Pi_1$ by an ortho\-gonally equivalent one, we may assume that $P_0 = R(\theta_1, \ldots, \theta_n)$ for some real numbers $-\pi \leqq \theta_1, \ldots, \theta_n < \pi$.

We shall presently prove the existence of an extension $\Pi: G \to \OG(2n)$ under the assumption that $-\pi < \theta_1, \ldots, \theta_n < \pi$. This assumption is certainly true when $q$ is odd, because of the identity ${P_0}^q = I_{2n}$. By the previously noticed surjectivity of the map $G_1 \subset G \onto G/G^{(e)}$, it will be enough to show that there is a continuous group homomorphism $\Upsilon: G^{(e)} \to \SO(2n)$ which restricts to $\Pi_1$ on $C = G^{(e)} \cap G_1 \subset G^{(e)}$ and satisfies the equation
\begin{equation}\label{equ:Upsilon-Pi1}
\Upsilon \bigl( g_1 x {g_1}^{-1} \bigr) = \Pi_1(g_1) \Upsilon(x) \Pi_1(g_1)^{-1} \quad \text{for all $x \in G^{(e)}$ and all $g_1 \in G_1$,}
\end{equation}
for then we can define $\Pi(g) := \Pi_1(g_1)\Upsilon(x)$ for all $g = g_1x \in G$. In order to construct $\Upsilon$, we observe that by our assumption on $\theta_1, \ldots, \theta_n$ and by Proposition \ref{prop:INJR} in combination with Remark \ref{rmk:X(theta1,...,thetan)} (see the appendix below) there exists a {\em unique} one-parameter subgroup
\begin{equation}\label{alpha}
\alpha: \nR \to \SO(2n) \quad \text{with} \quad \alpha(1) = R(\theta_1, \dotsc, \theta_n) = P_0 \quad \text{and} \quad \supnorm{\dot\alpha(0)} < \pi\text,
\end{equation}
namely, $t \mapsto \exp \bigl[ t \, X(\theta_1, \dotsc, \theta_n) \bigr]$ (the notations of the appendix are in use). As observed above, for each $g_1 \in G_1$ we must have $g_1x{g_1}^{-1} = x^{\pm 1}$ for all $x \in G^{(e)}$, the sign $\pm$ depending only on $g_1$. Accordingly, $\alpha(1)^{\pm 1} = {\Pi_1(g_1)\, \alpha(1)\, \Pi_1(g_1)^{-1}}$. Now, for each orthogonal matrix $R \in \OG(2n)$ the curves $t \mapsto {R\, \alpha(t) R^{-1}}$ and $t \mapsto \alpha(t)^{-1}$ define one-parameter subgroups ${R \alpha R^{-1}}$ and $\alpha^{-1}$ in $\SO(2n)$ for which the estimates $\supnorm{({R \alpha R^{-1}})'(0)} = \supnorm{\mathrm{Ad}(R) \cdot \alpha'(0)} \leqq \supnorm{\alpha'(0)} < \pi$ and $\supnorm{(\alpha^{-1})'(0)} = \supnorm{-\alpha'(0)} < \pi$ still hold, by Remark \ref{rmk:|Ad(R)(X)|}. By the uniqueness argument mentioned before, we see that $\alpha(t)^{\pm 1} = {\Pi_1(g_1)\, \alpha(t)\, \Pi_1(g_1)^{-1}}$ for all $t \in \nR$.
\end{paragrafo}

\begin{paragrafo}[The even-order case]\label{par:q=2p}
Suppose $q = \norm{C} = 2p$ for some $1 \leqq p \in \nN$. Define $N$ and $n$ as before. For each $g_1 \in G_1$ define the following subspace of $R = \C^0(G_1)$
\begin{equation}\label{def:Theta_g}
\Theta_{g_1} := \biggl\{ f \in \C^0(G_1) \biggl| \textsum[q]{i=1}{f(g_1 {c_0}^i)} = 0 \text{,~} \supp{f} \subset g_1C \biggr\}\text,
\end{equation}
$c_0$ being the chosen generator of the cyclic group $C$. Since $g_1C \subset H$ when $g_1 \in H$, and $g_1C \subset H'$ when $g_1 \in H'$, $\Theta_{g_1}$ is always a subspace of both $\Lambda_1$ and $\Lambda_2$ and hence so is
\begin{equation}\label{def:Theta}
\Theta := \bigoplus_{\widetilde{g_1} \in G_1/C} \Theta_{g_1} \subset R\text.
\end{equation}
The dimension of $\Theta$ is $[G_1:C](q-1) = (N/q)(q-1)$. The orthogonal complement of $\Theta_i$ in $\Lambda_i$ (the notation $\Theta_i$ simply indicates that we regard $\Theta$ as a subspace of $\Lambda_i$) is the subspace ${\Theta_i}' = \bigl\{ f \in \Lambda_i \bigl| f(g_1) = f(g_1c_0)\: \forall g_1 \in G_1 \bigr\}$ (of dimension $N/q-1$). Now,
\begin{elenco}
\item the generator $c_0 \in C$ fixes each $\Theta_{g_1}$, hence in particular the whole $\Theta$, and acts as the identity on the complement ${\Theta_i}'$; moreover,
\item for each $g, g_1 \in G_1$ and every $f \in \C^0(G_1)$ one has $f \in \Theta_g$ $\aeq$ ${g_1 \cdot f} \in \Theta_{g {g_1}^{-1}}$.
\end{elenco}
The second property follows by (\ref{equ:action-by-conj}) from the fact that ${g_1\, {c_0}^i\, {g_1}^{-1}} = {c_0}^{\pm i}$, the sign being a plus or a minus according to whether $g_1 \in H$ or $g_1 \in H'$.

Next let us concentrate on the action of $c_0$ on a single subspace $\Theta_g$ for $g \in G_1$ fixed. The $(-1)$-eigenspace for the action of $c_0$ on $\Theta_g$, namely
\begin{equation}\label{def:lambda_g}
\Lambda_g := \langle \lambda_g \rangle := \Bigl\langle \e{g} - \e{g c_0} + \cdots + \e{g {c_0}^{q-2}} - \e{g {c_0}^{q-1}} \Bigr\rangle = \biggl\langle \textsum[q]{i=1}{(-1)^i \e{g {c_0}^i}} \biggr\rangle \subset \Theta_g\text,
\end{equation}
is one-dimensional. Note that since $q = 2p$ is even, the $(-1)$-eigenvector $\lambda_g \in \Theta_g$ possesses the ``axial symmetry'' $\lambda_g(g {c_0}^i) = \lambda_g(g {c_0}^{-i})$. Thus ${g_1 \cdot \lambda_g} = \lambda_{g {g_1}^{-1}}$ for all $g, g_1 \in G_1$. Write down the orthogonal decomposition $\Theta_g = \Lambda_g \oplus {\Lambda_g}'$. It then follows that $g_1 \in G_1$ maps each complement ${\Lambda_g}'$ bijectively onto $(\Lambda_{g {g_1}^{-1}})'$. Now, introduce the following vectors
\begin{equation}\label{def:lambda_i^g}
\begin{cases}
\lambda_1^g := \lambda_g & \text{for each $g \in G_1$}
\\
\lambda_2^h := \lambda_h & \text{for each $h \in H$}
\\
\lambda_2^{h'} := -\lambda_{h'} & \text{for each $h' \in H'$.}
\end{cases}
\end{equation}
From the identity ${g_1 \cdot \lambda_g} = \lambda_{g {g_1}^{-1}}$ above, and the fact that ${g {g_1}^{-1}} \in H$ if and only if $g$ and $g_1$ both belong to $H$ or both belong to $H'$, we obtain the following transformation rules: (a)~${g_1 \cdot \lambda_1^g} = \lambda_1^{g {g_1}^{-1}}$ for all $g, g_1 \in G_1$; (b)~${h \cdot \lambda_2^g} = \lambda_2^{g h^{-1}}$ for all $h \in H$, $g \in G_1$; (b$'$)~${h' \cdot \lambda_2^g} = -\lambda_2^{g {h'}^{-1}}$ for all $h' \in H'$, $g \in G_1$. Put $\Lambda_i^g := \bigl\langle \lambda_i^g \bigr\rangle$ ($i = 1,2$) and
\begin{equation}\label{def:Lambda^g}
\Lambda^g := \Lambda_1^g \oplus \Lambda_2^g \subset \Lambda\text.
\end{equation}
On this two-dimensional subspace of $\Lambda$ we fix the basis $(\lambda_1^g,0)$, $(0,\lambda_2^g)$.

We have the following decomposition of $\Lambda$ into three $G_1$-invariant direct summands
\begin{equation}\label{equ:split-Lambda}
\Lambda = \Biggl[\; \bigoplus_{\widetilde{g_1} \in G_1/C} \Lambda^{g_1}\; \Biggr] \oplus \Biggl[\; \bigoplus_{\widetilde{g_1} \in G_1/C} \bigl( {\Lambda_{g_1}}' \oplus {\Lambda_{g_1}}' \bigr)\; \Biggr] \oplus \Bigl[\Theta_1' \oplus \Theta_2'\Bigr]\text.
\end{equation}
As in the preceding subsection \ref{par:q=2p+1}, we are given a representation $\Pi_1: G_1 \to \GL(\Lambda)$ (namely, the $G_1$-module $\Lambda$) which we would like to extend to a continuous representation $\Pi: G \to \GL(\Lambda)$ of the whole $G$ on the same vector space. Of course, to this end, we may deal with each direct summand in (\ref{equ:split-Lambda}) separately.

The right-hand summand in (\ref{equ:split-Lambda}) is a trivial $G_1$-module, so on that summand we may extend $\Pi_1$ by the trivial representation.

The mid summand falls into the case already studied in Subsection \ref{par:q=2p+1} because it corresponds to an orthogonal $G_1$-action in which $c_0$ acts with no $-1$ eigenvalue. By reasoning as before one can obtain an extension on that summand by using the results of Appendix \ref{app:INJR}.

As to the first summand, let us call it $W$ for brevity, we construct an extension directly. Define a one-parameter subgroup $\alpha: \nR \to \GL(W)$ as follows:
\begin{equation}\label{def:alpha:-1}
\alpha(t) := \bigoplus_{\widetilde{g_1} \in G_1/C} \alpha^{g_1}(t)
\end{equation}
where $\alpha^{g_1}(t) \in \GL(\Lambda^{g_1})$ is the linear map represented by the $2 \times 2$ matrix $\bigl( \begin{smallmatrix} \cos t & -\sin t \\ \sin t & \phantom-\cos t \end{smallmatrix} \bigr)$ with respect to the preferred basis $(\lambda_1^{g_1},0)$, $(0,\lambda_2^{g_1})$ of $\Lambda^{g_1}$. Note that the linear map $\alpha^{g_1}(t)$ is well-defined, i.e., independent of the choice of a representative $g_1$ for the given coset $\widetilde{g_1} \in G_1/C$. Indeed, $\lambda_i^{g_1c_0} = -\lambda_i^{g_1}$ for $i=1,2$ so that the preferred basis associated with $g_1c_0$ is minus the one associated with $g_1$ and therefore the same matrix represents the same linear map. Observe that $\alpha(\pi) = -\id = \Pi_1(c_0)|_W$. As in Subsection \ref{par:q=2p+1}, we need to check that
\begin{equation}\label{alpha-pi}
\begin{cases}
\alpha(t) = \Pi_1(h)\, \alpha(t)\, \Pi_1(h)^{-1} & \text{for all $h \in H$}
\\
\alpha(t)^{-1} = \Pi_1(h')\, \alpha(t)\, \Pi_1(h')^{-1} & \text{for all $h' \in H'$.}
\end{cases}
\end{equation}
For each $g_1, g \in G_1$ the linear automorphism $\Pi_1(g)|_W \in \GL(W)$ maps the subspace $\Lambda^{g_1}$ onto the subspace $\Lambda^{g_1 g^{-1}}$. On the other hand $\alpha(t)$ maps $\Lambda^{g_1}$ and $\Lambda^{g_1 g^{-1}}$ into themselves by construction. Thus we can check the identities (\ref{alpha-pi}) for $g=h,h'$ in the preferred bases of $\Lambda^{g_1}$ and $\Lambda^{g_1 g^{-1}}$. This is straightforward, in view of the transformation rules (a), (b) and (b$'$) stated immediately after (\ref{def:lambda_i^g}).
\end{paragrafo}

\appendix

\section[Exponential injectivity radius]{Computation of the exponential injectivity radius for the classical orthogonal linear groups}\label{app:INJR}
The Lie group $G = \SO(n)$ is semisimple for all $n \geqq 3$. In fact, $\mathit Z\bigl( \SO(2p-1) \bigr) = \{I\}$ and $\mathit Z\bigl( \SO(2p) \bigr) = \{\pm I\}$ for all $p \geqq 2$. Compare \cite[Remark (3.14) p.~201]{BtD95}, and \cite[p.~102]{FuHa91}. Thus \mbox{Lazard}--{Tits'} criterion \cite[Th\'eor\`eme (2.1)]{LazTits66}, which requires the connected component of the centre $\mathit Z_0(G)$ to be simply connected, can be applied when $G = \SO(n)$.

The Lie algebra $\mathfrak g = \so(n)$ is given by the $n \times n$ skew-symmetric real matrices $X \in \mathrm{Mat}_{n \times n}(\nR)$, $\transposed X = -X$, the Lie bracket being the usual anticommutator $[X,Y] = XY-YX$. Now make the identification $\mathrm{Mat}_{n \times n}(\nR) = \End(\nR^n)$, put the euclidean norm on $\nR^n$, and let \supnorm{X} denote the usual operator norm on $\End(\nR^n)$. Since for this norm one has the inequality $\supnorm{XY} \leqq {\supnorm X\supnorm Y}$, one obtains the following estimate for the Lie bracket
\begin{equation}\label{est:||[X,Y]||}
\bigl\|[X,Y]\bigr\| \leqq 2\supnorm{X}\supnorm{Y}\text.
\end{equation}
We will see that the factor two here is actually the norm of the bilinear form $(X,Y) \mapsto [X,Y]$, in other words (\ref{est:||[X,Y]||}) is the best estimate possible. If we put $\norm{X} := 2\supnorm{X}$, the inequality (\ref{est:||[X,Y]||}) implies that the latter norm (on the Lie algebra $\mathfrak g$) is \textit{admissible,} i.e.\ satisfies
\begin{equation}\label{est:|[X,Y]|}
\bigl|[X,Y]\bigr| \leqq \norm{X}\norm{Y}\text.
\end{equation}
It then follows from \cite[Th\'eor\`eme (2.1)]{LazTits66} that the exponential mapping $X \mapsto \exp(X) = I + X + X^2/2 + \cdots$ from the Lie algebra $\mathfrak g = \so(n)$ into the respective Lie group $G = \SO(n)$ is injective on the open ball $\{X \in \mathfrak g: \norm{X} < \pi\} = \{X \in \mathfrak g: \supnorm{X} < \pi/2\}$. It is the goal of the present appendix to show that the injectivity radius of the exponential mapping $\exp: \so(n) \to \SO(n)$ is actually twice as much:

\begin{proposizione}\label{prop:INJR}
The injectivity radius of the exponential mapping, with respect to the admissible norm $\norm{X} = 2\supnorm{X}$ on the Lie algebra $\mathfrak g = \so(n)$ of the special orthogonal group $G = \SO(n)$, is exactly $2\pi$.
\end{proposizione}
Equivalently, the injectivity radius of the exponential mapping with respect to the operator norm \supnorm{X} on $\mathfrak g$ is exactly $\pi$. This is actually the best we can hope for, in view of the first of the following two remarks (which we will need in the proof of the Proposition):

\begin{paragrafo}[Remark]\label{rmk:X(theta1,...,thetan)}
The skew-symmetric $2 \times 2$ matrix $X(\theta) = \bigl( \begin{smallmatrix} 0 & -\theta \\ \theta & \phantom-0 \end{smallmatrix} \bigr)$ exponentiates to $R(\theta) = \bigl( \begin{smallmatrix} \cos\theta & -\sin\theta \\ \sin\theta & \phantom-\cos\theta \end{smallmatrix} \bigr)$. More generally, the skew-symmetric $2n \times 2n$ matrix
$$%
X(\theta_1, \ldots, \theta_n) = \begin{pmatrix} X(\theta_1) & \cdots & 0 \\ \vdots & \ddots & \vdots \\ 0 & \cdots & X(\theta_n) \end{pmatrix} \quad \text{exponentiates to} \quad \begin{pmatrix} R(\theta_1) & \cdots & 0 \\ \vdots & \ddots & \vdots \\ 0 & \cdots & R(\theta_n) \end{pmatrix}\text.
$$%
Note that $\supnorm{X(\theta_1, \ldots, \theta_n)} = {\underset{i = 1, \ldots, n}{\max} \norm{\theta_i}}$.
\end{paragrafo}

\begin{paragrafo}[Remark]\label{rmk:|Ad(R)(X)|}
The orthogonal matrices $R \in \OG(n)$, $R \cdot \transposed R = I$ act on the Lie algebra $\so(n)$ by $X \mapsto RXR^{-1}$. In fact, $RXR^{-1} = {\mathrm{Ad}(R) \cdot X}$ is precisely the adjoint representation $\mathrm{Ad}: \OG(n) \to \GL\bigl(\so(n)\bigr)$. Notice that each orthogonal matrix $R \in \OG(n)$ preserves the euclidean distance as an operator $\in \End(\nR^n)$ and therefore $\supnorm{R} = \supnorm{\transposed R} = 1$. It follows that
\begin{equation}\label{equ:||Ad(R)(X)||}
\supnorm{\mathrm{Ad}(R) \cdot X} = \supnorm{R \cdot X \cdot \transposed R} \leqq \supnorm{R}^2\supnorm{X} = \supnorm{X}\text.
\end{equation}
Hence for each $R \in \OG(n)$ the adjoint automorphism $\mathrm{Ad}(R)$ on $\mathfrak g = \so(n)$ maps the open ball $\mathit B(0,r) \subset \mathfrak g$ of radius $r > 0$ into itself. The same remark applies of course to the other norm: $\norm{\mathrm{Ad}(R)(X)} \leqq \norm{X}$ for all $X \in \mathfrak g$, $R \in \OG(n)$.
\end{paragrafo}

\begin{proof}[Proof of Proposition \ref{prop:INJR}]
For simplicity we shall assume that is $n$ even, as this is the only case of practical interest for us. (There are obvious modifications for $n$ odd.) So, let $X, Y \in \mathfrak g = \so(2n)$ be given with $\norm{X}, \norm{Y} < 2\pi$, and suppose $\exp(X) = \exp(Y) \in \SO(2n)$.

Since the norm \norm{X} on $\mathfrak g$ is admissible, it follows from \cite[Lemme (6.1)]{LazTits66} that $[X,Y] = 0$. Therefore, the linear subspace $\mathrm{Span}\{X,Y\}$ is an abelian subalgebra of $\mathfrak g$ and hence there is a maximal abelian subalgebra $\mathfrak t \subset \mathfrak g$ such that $X, Y \in \mathfrak t$. Let $T = \exp(\mathfrak t)$ be the maximal torus in $\SO(2n)$ corresponding to this maximal abelian subalgebra. The one-parameter subgroups $\alpha_X := \bigl\{ t \mapsto \exp(tX) \bigr\}$ and $\alpha_Y$ are contained in the maximal torus $T$. Since all maximal tori in a connected compact Lie group are conjugated to each other \cite[Theorem (1.6) p.~159, and p.~5]{BtD95}, there will be an element $g_0 \in \SO(2n)$ with ${g_0 T {g_0}^{-1}} = T(n)$, where $T(n)$ denotes the standard torus in $\SO(2n)$ consisting of all block-diagonal matrices of the form $R(\theta_1, \ldots, \theta_n)$. Since by Remark \ref{rmk:|Ad(R)(X)|} the linear automorphism $\mathrm{Ad}(g_0) \in \GL(\mathfrak g)$ maps the open ball $\mathit B(0,2\pi)$ into itself, if we put $X_0 = \mathrm{Ad}(g_0)(X)$ and $Y_0 = \mathrm{Ad}(g_0)(Y)$ then we still have $\norm{X_0}, \norm{Y_0} < 2\pi$. Clearly, $X=Y$ if and only if $X_0=Y_0$.

Now, the one-parameter subgroups $\alpha_{X_0} = {g_0 \alpha_X {g_0}^{-1}}$ and $\alpha_{Y_0} = {g_0 \alpha_Y {g_0}^{-1}}$ (recall that by the naturality of the exponential mapping one has
$$%
\exp\bigl( \mathrm{Ad}(g_0) \cdot X \bigr) = {g_0\, \exp X\, {g_0}^{-1}}\text,
$$%
compare \cite[(3.2) p.~23]{BtD95} for instance) are contained in the torus $T(n)$, and $\supnorm{X_0}, \supnorm{Y_0} < \pi$. From Remark \ref{rmk:X(theta1,...,thetan)} it follows that $X_0 = Y_0$ and hence that $X=Y$.
\end{proof}

\hbadness=10000

\end{document}